\documentclass[10pt]{article}

\usepackage[dvips]{graphicx}
\usepackage{amssymb,amsfonts,amsmath,amsthm}
\usepackage{psfrag}
\usepackage[numbers]{natbib}

\newcommand{\power}{\mathscr{P}}
\newcommand{\fs}[1]{\mathcal{#1}}
\newcommand{\nats}{\mathbb{N}}
\newcommand{\reals}{\mathbb{R}}
\newcommand{\order}{\mathcal{O}}
\newcommand{\define}{:=}
\newcommand{\comment}[1]{}
\newcommand{\editor}[1]{\texttt{#1}}
\newcommand{\expect}{\mathbf{E}}
\newcommand{\var}{\mathbf{VAR}}

\theoremstyle{plain}
\newtheorem{theorem}{Theorem}

\newtheorem{corollary}[theorem]{Corollary}

\theoremstyle{definition}
\newtheorem{definition}{Definition}

\newcommand\independent{\protect\mathpalette{\protect\independenT}{\perp}} \def\independenT#1#2{\mathrel{\rlap{$#1#2$}\mkern2mu{#1#2}}}

\begin{document}

%%%%%%%%%%%%%%%%%%%%%%%%%%%%%%%%%%%%%%%%%%%%%%%%%%%%%%%%%%%%%%%%%%%%
% Header                                                           %
%%%%%%%%%%%%%%%%%%%%%%%%%%%%%%%%%%%%%%%%%%%%%%%%%%%%%%%%%%%%%%%%%%%%

\title{Thermodynamics as a theory of decision-making with information processing costs}
\author{Pedro A. Ortega and Daniel A. Braun}

\maketitle

\begin{abstract}
Perfectly rational decision-makers maximize expected utility,
but crucially ignore the resource costs incurred when determining
optimal actions. Here we propose an information-theoretic formalization of
bounded rational decision-making where decision-makers trade off
expected utility and information processing costs. Such bounded rational
decision-makers can be thought of as thermodynamic machines that undergo
physical state changes when they compute. Their behavior is governed by a free energy functional that trades off changes in internal energy---as a proxy for utility---and
entropic changes representing computational costs induced by changing states. As a result,
the bounded rational decision-making problem can be rephrased in terms of well-known concepts
from statistical physics.
In the limit when computational costs are ignored, the maximum
expected utility principle is recovered.
We discuss the relation to \emph{satisficing} decision-making procedures as well as links to existing theoretical frameworks and human decision-making experiments that describe deviations from expected utility theory.
Since most of the mathematical machinery can be borrowed from statistical physics,
the main contribution is to axiomatically derive and interpret the thermodynamic free energy
as a model of bounded rational decision-making.
\end{abstract}

%%%%%%%%%%%%%%%%%%%%%%%%%%%%%%%%%%%%%%%%%%%%%%%%%%%%%%%%%%%%%%%%%%%%
% Body                                                             %
%%%%%%%%%%%%%%%%%%%%%%%%%%%%%%%%%%%%%%%%%%%%%%%%%%%%%%%%%%%%%%%%%%%%

\section{Introduction}
In everyday life decision-makers often have to make fast and frugal choices \cite{Gigerenzer1999,Gladwell2005}. Consider, for example, an antelope that quickly has to choose a direction of flight when faced with a predator. By the time an antelope had considered all possible flight paths to determine the optimal one, it would most probably be already eaten. In general, decision-makers seem to trade off the expected desirability of the consequences of an action against the effort and resources (time, money, food, computational effort, knowledge, opportunity costs, etc.) required for searching the optimum \cite{NivDawJoelEtAl07,daw2012}.

Classic theories of decision making generally ignore information-processing costs by assuming that decision makers always pick the option with maximum return---irrespective of the effort or the resources it might take to find or compute the optimal action \cite{Neumann1944,Savage1954,Fishburn1982}. Such decision-makers are described as \emph{perfectly rational}. However, being perfectly rational seems to contradict our intuition of real-world decision-making, where information processing constraints play an important role \cite{Gigerenzer1999}.
This has led to an abundant literature on \emph{bounded rationality} \cite{Simon1972,Simon1984,Rubinstein1998,Gigerenzer2001}. Unlike perfectly rational decision makers, bounded rational decision-makers are subject to information processing constraints, that is they may have limited time and speed to process a limited amount of information.

\subsection{Thermodynamic Intuition}

\begin{figure}[htbp]
\begin{center}
    \small
    \psfrag{l1}[c]{a)}
    \psfrag{l2}[c]{b)}
    \psfrag{c1}[l]{$p V$}
    \psfrag{c2}[l]{$(1-p) V$}
    \psfrag{p1}[c]{A}
    \psfrag{p2}[c]{B}
    \includegraphics[]{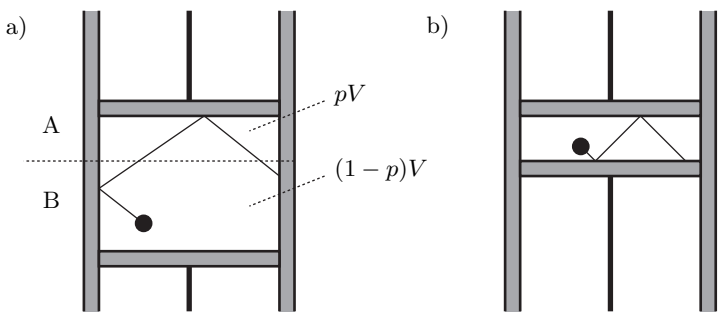}
    \caption[The Molecule-In-A-Box Device.]{The Molecule-In-A-Box Device.
    (a) Initially, the molecule moves freely within a space of volume $V$ delimited by two pistons. The compartments A and B correspond to the two logical states of the device. (b) Then, the lower piston pushes the molecule into part A having volume $V' = p V$.}
    \label{fig:molecule-in-a-box}
\end{center}
\end{figure}

Here we follow a thermodynamic
argument \cite{Feynman1996} that allows measuring resource (or information) costs in physical
systems in units of energy. The generality of the argument relies on the fact that ultimately any real agent has to be incarnated in a physical system, as any process of information processing must always be accompanied by a pertinent physical process \cite{Tribus1971}.
In the following we conceive of information processing as changes in information states (i.e. ultimately changes of probability distributions),
which consequently implies changes in physical states, such as flipping gates in a transistor, changing voltage on a microchip, or even changing location of a gas particle.
Such state changes in physical systems are not for free, that is the do not happen spontaneously.
Consequently, if we want to control a physical system into a desirable state we also have
to take into consideration that changing from the current state to the desirable state
incurs a cost.

According to Landauer's principle, one can postulate a formal correspondence between one unit of information and one unit of energy \cite{Landauer1961,Bennett1973,Bennett1982}. Consider representing one bit of information using one of the following logical devices: a molecule that can be located either on the top or the bottom part of a box; a coin whose face-up side can be either head or tail; a door that can be either open or closed; a train that can be orientated facing either north or south; and so forth. Assume that all these devices are initialized in an undetermined logical state, where the first state has probability $p$ and the second probability $1-p$. Now, imagine you want to set these devices to their first logical state. In the case of the molecule in a box, this means the following. Initially, the molecule is uniformly moving around within a space confined by two pistons as depicted in Figure~\ref{fig:molecule-in-a-box}a. Assuming that the initial volume is $V$, the molecule has to be pushed by the lower piston into the upper part of the box having volume $V' = p V$ (Figure~\ref{fig:molecule-in-a-box}b). From information theory, we know that the number of bits that we fix by this operation is given by $- \log p$.

To make things concrete, we assume that the device has diathermal walls and is immersed in a heat bath at constant temperature $T$. Since the walls are diathermal, the temperature inside of the box is maintained at the temperature of the heat bath. We model the particle as an ideal gas. When an ideal gas is compressed under isothermal conditions from an initial volume $V$ to a final volume $V'$, then the work is calculated as
\begin{equation}\label{eq:mechanical-work}
    W
    = - \int_{V}^{V'} \frac{ NkT }{ V } \, dV
    = NkT \ln \frac{ V }{ V' },
\end{equation}
where $N \geq 0$ is the amount of substance and $k > 0$ is the Boltzmann constant. The minus sign is just a convention to denote work done by the piston rather than by the gas. If we assume $N=1$ and make use of the fact that $V'=pV$ we get
\[
    W
    = kT \ln \frac{ V }{ p V }
    = -kT \ln p = -\frac{kT}{\log e} \log p
    = -\gamma_\text{mol} \log p,
\]
where the constant $\gamma_\text{mol} \define
\frac{RT}{\log e} > 0$ can be interpreted as the conversion factor between one
unit of information and one unit of energy for the molecule-in-a-box device.

How do we compute the information and work for the case of the coin, door and
train devices? The important observation is that we can model these cases as if
they were like molecule-in-a-box devices, with the difference that
their conversion factors between units of information and units of work are
different. Hence, the number of bits fixed while these devices are set to the
first state is given by $- \log p$, i.e.\ exactly as in the case of the molecule. However, the work is given by
\[
    -\gamma_\text{coin} \log p,
    \quad -\gamma_\text{door} \log p,
    \quad \text{and}
    \quad -\gamma_\text{train} \log p
\]
respectively, where $\gamma_\text{coin}$, $\gamma_\text{door}$ and
$\gamma_\text{train}$ are the associated conversion factors between units of information. Obviously,
$
    \gamma_\text{mol}
    \leq \gamma_\text{coin}
    \leq \gamma_\text{door}
    \leq \gamma_\text{train}.
$
The point is that changes in knowledge states are costly and that these costs are proportional to the information. In the next section, we derive
a general expression of information costs in physical systems that make decisions.

\section{Information-Theoretic Foundations} \label{sec:foundation}

\subsection{Resource Costs}

We model any observable sequential process, such as a sequence of
interactions or a sequence of computation steps, as a filtration on a
measure space. To simplify our exposition, we consider only finite
measure spaces. Let $(\Omega,\Sigma)$ denote a measurable space, where
$\Omega$ denotes the sample space and where $\Sigma$ is a
$\sigma$-algebra on $\Omega$. Let $p$ be a conditional probability
measure on $(\Omega,\Sigma)$, such that for any two events $A, B \in
\Sigma$, $p(A|B)$ denotes the conditional probability of the~$A$
given~$B$, where the condition~$B$ plays the role of the current
information state of the process. The sequential realization of a
process is modelled as a sequence of conditions $A_1, A_2, \ldots,
A_T$ on the sample space $\Omega$, where each new condition $A_t$
refines the current information state $\bigcap_{\tau \leq t} A_\tau$
by excluding the complement~$A_t^\complement$.

We further assume that a transformation of an information state from
$B$ to $(A \cap B)$ entails a cost $\rho(A|B)$ that could be measured
in dollars, time or any arbitrary scale of effort. Moreover, we assume
that this transformation cost is decomposable; that is, if we undergo
a knowledge change from $C$ to $(A \cap B \cap C)$, then we should pay
the same cost as undergoing a change first from $C$ to $(B \cap C)$
and then from $(B \cap C)$ to $(A \cap B \cap C)$. Finally, the
quintessential  information-theoretic postulate is that conditional
probabilities impose a monotonic order over transformation
costs\footnote{This intuition is central for optimal coding theory
where short codewords are assigned to frequent events and long
codewords are assigned to rare events \cite{Mackay2003}. Therefore, we
could regard the codeword length as a valuable resource that we have
to bet on events with different probabilities.}. We can sum up our
postulates as follows:

\begin{definition}[Axioms of Transformation Costs]\label{def:cost}
Let $(\Omega, \Sigma)$ be a measurable space and let
$p:(\Sigma\times\Sigma) \rightarrow [0,1]$ be a conditional
probability measure over $\Sigma$ (i.e. for any $A \in \Sigma$,
$p(\cdot|A)$ is a probability measure over $A$). A function $\rho:
(\Sigma \times \Sigma) \rightarrow \reals^+$ is a transformation cost
function for $p$ iff it has the following three properties for all
events $A, B, C, D \in \Sigma$:
\begin{align*}
    &\text{A1.~real-valued:}
        &&\exists f, \quad \rho(A|B) = f\bigr( p(A|B) \bigr)\in \reals,\\
    &\text{A2.~additive:}
    &&\rho(A \cap B|C) = \rho(A|C) + \rho(B|A \cap C),\\
    &\text{A3.~monotonic:}
    &&[\rho(A|B) > \rho(C|D)] \quad \Leftrightarrow \quad [p(A|B)
\lessgtr p(C|D)].
\end{align*}
\end{definition}

These three properties enforce a strict correspondence between
probabilities and transformation costs
\cite{OrtegaBraun2010b,Ortega2011}.

\begin{theorem}[Transformation Costs $\leftrightarrow$
Probabilities]\label{theo:utility}
If $f$ is such that $\rho(A|B) = f(p(A|B))$ for every choice of the
probability space $(\Omega, \Sigma, p)$, then $f$ is of the form
\[
    f(\cdot) = -\tfrac{1}{\beta} \log(\cdot),
\]
where $\beta$ is a real parameter.
\end{theorem}

That is, the transformation cost $\rho(A|B)$ is \emph{proportional} to
the information content $-\log p(A|B)$, where the parameter $\beta$
plays the role of the conversion factor. The logarithmic mapping
between probabilities and ``costs'' is well-known in information
theory, and there are many possible ways to derive it
\cite{Shannon1948,Csiszar2008}. The important observation is that our
derivation stems purely from postulates regarding transformation
costs.

According to Definition~\ref{def:cost}, transformation costs measure
the relative cost of an event \emph{relative} to a reference event.
However, we can also introduce an \emph{absolute} cost measure to
single events such that transformation costs are obtained as
differences.

\begin{definition}[Potential]
\label{def:potential}
Let $\rho$ be a transformation cost function. A set function
$\phi:\Sigma \rightarrow \reals$ is called a \emph{cost potential} for
$\rho$ iff for all $A, B \in \Sigma$,
\begin{align*}
        \phi(\Omega) &:= \phi_0\\
        \phi(A \cap B) &:= \phi(B) + \rho(A|B) \qquad \forall A, B \in \Sigma,
\end{align*}
where $\phi_0$ is an arbitrary real value.
\end{definition}
One can easily verify that this potential is well defined for all
events, and that $\rho(A|B) = \phi(A \cap B) - \phi(B)$. It captures
the intuition that starting out from the high-probability event $B$
with potential $\phi(B)$ one has to pay the cost $\rho(A|B)$ to arrive
at the low-probability event $A \cap B$ with potential $\phi(A \cap
B)$.

In the following, consider a reference set $S \in \Sigma$ having a
measurable partition~$\fs{X}$. Cost potentials have an important
recursive structure: the cost potential of an event is uniquely
determined by the potential of its constituent events. If~$\fs{X}$ is
a measurable partition of a reference event $S \in \Sigma$, then
\begin{equation}\label{eq:log-partition}
        \phi(S) = -\tfrac{1}{\beta} \log\sum_{x \in \fs{X}} e^{-\beta \phi(x)}.
\end{equation}
Furthermore, the probability of a member $x \in \fs{X}$ of the
partition relative to $S$ can be expressed as a \emph{Gibbs measure}:
\begin{equation}\label{eq:gibbs}
        p(x|S)
    = \frac{e^{-\beta \phi(x)}}{e^{-\beta \phi(S)}}
    = \frac{e^{-\beta \phi(x)}}{\sum_{x \in \fs{X}} e^{-\beta \phi(x)}}.
\end{equation}
 In statistical physics it is well-known that the Gibbs measure
satisfies a variational principle in the \emph{free energy}, which is
defined as
\begin{equation}\label{eq:free-energy}
    F_\beta[q] := \sum_{x \in \fs{X}} q(x) \phi(x)
        + \frac{1}{\beta} \sum_{x \in \fs{X}} q(x) \log q(x).
\end{equation}
More specifically, it is well known that for any probability measure
$q$ over the partition $\fs{X}$ of $S$,
\begin{equation}\label{eq:variational-principle}
    F[q] \geq F[p]
    = -\frac{1}{\beta} \log \phi(S),
\end{equation}
where the lower bound is attained by the Gibbs measure $p(x) \propto
e^{-\beta \phi(x)}$. Equations~\eqref{eq:log-partition}
to~\eqref{eq:variational-principle} constitute fundamental results
that will be generalized and interpreted in the next section.

\subsection{Gains and Losses}

Equipped with the results from the preceding section, we can now
proceed to model a bounded rational decision maker. Because
transformation costs matter, we model a decision as a transformation
of a prior behavior into a final behavior, where we represent the
direction of change as a utility criterion.

The Gibbs measure in~\eqref{eq:gibbs} allows us describing a
probability measure $p$ over a partition $\fs{X}$ in terms of a cost
potential $\phi$ over $\fs{X}$. In particular, we see that a
decision-maker's a priori behavior or belief described by $p_0(x)$ and
$\phi_0(x)$ changes to $p(x)$ and $\phi(x)$ if he is exposed to the
gain (or loss) $U(x)$, such that
\begin{equation}
\phi(x) = \phi_0(x)-U(x)
\end{equation}
and
\begin{equation}
\label{eq:probdist}
p(x) \propto e^{-\beta \phi_0(x) + \beta U(x) } \propto p_0(x) e^{\beta U(x)}
\end{equation}
as illustrated in Figure~1. The function $U$ represents either gains
or losses and not absolute levels of costs, because it expresses a
difference in the potential $U(x)=\phi_0(x)-\phi(x)$. The equilibrium
distribution~\eqref{eq:probdist} that arises in a change can also be
characterized in terms of a variational principle, in a manner
analogous to~\eqref{eq:variational-principle}.

\begin{theorem}[Negative Free Energy
Difference]\label{theo:free-energy-difference}
Let $p_0(x)$ and $p(x)$ be the Gibbs measures with potentials
$\phi_0(x)$ and $\phi(x)$ and resource parameter $\beta$. Let $F_0$
and $F$ be the free energies minimized by $p_0$ and $p$ respectively.
Then, the negative free energy difference $-\Delta F = F_0 - F$ is
\begin{equation}
\label{eq:free-energy-difference}
    -\Delta F
    = \sum_{x \in \fs{X}} p(x) U(x)
        - \frac{1}{\beta} \sum_{x \in \fs{X}} p(x) \log \frac{p(x)}{p_0(x)},
\end{equation}
where $U(x) = \phi_0(x)-\phi(x)$.
\end{theorem}

Since the difference in the negative free energy $-\Delta F = F - F_0$
has the same dependency on $p$ as the free energy $F$, we can use
$-\Delta F$ directly as a variational principle in $p$.

\begin{corollary}[Variational Principle]
\label{theo:var-principle}
The negative free energy difference provides a variational principle
for the equilibrium distribution, i.e.\
\[
    -\Delta F[q]
    := \sum_{x \in \fs{X}} q(x) U(x)
    - \frac{1}{\beta} \sum_{x \in \fs{X}} q(x) \log \frac{q(x)}{p_0(x)}
\]
is maximized by
\[
    p(x) = \frac{1}{Z} p_0(x) e^{\beta U(x)},
    \quad \text{where} \quad
    Z \define \sum_{x \in \fs{X}} e^{\beta U(x)}.
\]
Furthermore,
\[
    \Delta F[q]
    \leq \Delta  F[p]
    = \tfrac{1}{\beta} \log Z.
\]
\end{corollary}

\subsection{Choice \& Belief Probabilities}
The distribution \eqref{eq:probdist} can be interpreted both as an
action or observation probability in the context of bounded rational
decision-making. In the case of actions, $p_0$ represents the a priori
choice probability of the agent which is refined to the choice
probability $p$ when evaluating the imposed gain (or loss) $U$. The
associated change in probability depends on the resource parameter
$\beta$ and corresponds to the computation that is necessary to
evaluate the gains (or losses). In the case of observations, $p_0$
represents the a priori belief of the agent given by a probabilistic
model, which is then distorted due to the presence of possible gains
(or losses) that are evaluated by the holder of the belief. This way,
model uncertainty and risk-aversion can be parameterized by $\beta$.

\begin{figure}[tb]
\begin{center}
    \vspace*{.05in}
    \footnotesize
    \psfrag{a1}[c]{$\phi_0$}
    \psfrag{a2}[c]{$-U$}
    \psfrag{a3}[c]{$\phi = \phi_0 - U$}
    \psfrag{a4}[c]{Initial}
    \psfrag{a5}[c]{Final}
    \psfrag{a6}[l]{low}
    \psfrag{a7}[l]{high}
    \psfrag{a8}[c]{Probability}
    \includegraphics[width=8cm]{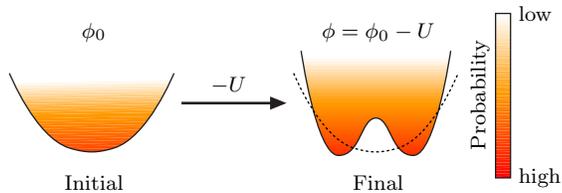}
    \caption{Representing a decision maker as a thermodynamic system,
the behavior of the decision-maker exposed to a gain $U$ can be
expressed as
    a change of his initial cost potential $\phi_0$ to a final cost
potential $\phi$, where $\phi = \phi_0-U$. The choice or belief
probabilities of the decision-maker change according to
\eqref{eq:probdist} from~$p_0$ to~$p$.}
    \label{fig:transformation}
\end{center}
\end{figure}

For different values of $\beta$ the distribution \eqref{eq:probdist}
has the following limits
\begin{eqnarray}
\lim_{\beta \rightarrow \infty} p(x) &=& \delta(x - x^*), \qquad x^* =
\max_x U(x) \nonumber \\
\lim_{\beta \rightarrow 0} p(x) &=& p_0(x) \nonumber \\
\lim_{\beta \rightarrow -\infty} p(x) &=& \delta(x - x^*), \qquad x^*
= \min_x U(x) \nonumber .
\end{eqnarray}
In the case of actions the three limits imply the following: The limit
$\beta \rightarrow \infty$ corresponds to the perfectly rational actor
that infallibly selects the action that maximizes gain (or minimizes
loss $-U(x)$.
The limit $\beta \rightarrow 0$ is an actor without resources that
simply selects his action according to his prior. The limit $\beta
\rightarrow -\infty$ corresponds to an actor that is
perfectly ``anti-rational'' and always selects the action with the
worst outcome. In the case of observations the three limits correspond
to an extremely optimistic observer ($\beta \rightarrow \infty$) who
believes only in the best possible outcome, an extremely pessimistic
observer ($\beta \rightarrow -\infty$) who anticipates only the worst,
and a risk-neutral Bayesian observer ($\beta \rightarrow 0$) who
simply relies on the probabilistic model $p_0$.

\subsection{The Certainty Equivalent}
In statistical physics \cite{Callen1985}, the free energy difference
\[
\Delta A = \Delta E - Q = W
\]
measures the amount of available ``good energy'' (work $W$) by
subtracting the ``bad energy'' (heat $Q$)
from the total energy $\Delta E = \expect[U]$.
%\footnote{Note that $A = -F$; i.e.\ we have chosen to state the
%results in terms of the negative free energy because the
%interpretation in economic terms is more natural.}.
The crucial physical intuition is that we have uncertainty about some
aspects of the objects that make up the heat energy, for example we do
not know the exact trajectories of all gas particles at temperature
$\beta$. This uncertainty means that we do not have full control over
the objects and cannot extract all the energy as work
\cite{Feynman1996}. Economically speaking, the physical concept of
work, and therefore also the difference in free energy, measures the
certainty equivalent of a gain (or loss) that is contaminated by
uncertainty. In general, we can therefore use the free energy
difference to ascribe a certainty equivalent value to choice
situations of the form \eqref{eq:probdist}. As can be seen from
Corollary~\ref{theo:var-principle}, this value is given by the log
partition function, i.e. the logarithm of the normalization constant
$Z$. For different values of $\beta$, the certainty equivalent takes
the following limits
\begin{eqnarray}
\lim_{\beta \rightarrow \infty} \frac{1}{\beta} \log Z &=& \max_x U(x)
\nonumber \\
\lim_{\beta \rightarrow 0} \frac{1}{\beta} \log Z &=& \sum_x p_0(x)
U(x) \nonumber \\
\lim_{\beta \rightarrow -\infty} \frac{1}{\beta} \log Z &=& \min_x
U(x) \nonumber .
\end{eqnarray}
Again, the case $\beta \rightarrow \infty$ corresponds to the
perfectly rational actor (or the extremely optimistic observer), the
case $\beta \rightarrow -\infty$ corresponds to the perfectly
``anti-rational'' actor (or the extremely pessimistic observer) and
the case $\beta \rightarrow 0$ corresponds to the actor that has no
resources (or the risk-neutral observer) such that the best one can
expect is the expected gain or loss.

Corollary~\ref{theo:var-principle} has two interpretations in
statistical physics, either as an instantiation of a  \emph{minimum
energy principle} or as a \emph{maximum entropy principle}
\cite{Callen1985}.
Accordingly, \eqref{eq:probdist} can either be seen as the
distribution that maximizes the entropy given a constraint on the
expectation value of $U$ or as the distribution that minimizes the
expectation of $-U$ given a constraint on the entropy of $p$. In the
context of observer modeling, the first interpretation provides a
principle for estimation and the second interpretation provides a
principle for bounded rational decision-making in the case of acting,
which is a maximum expected gain principle with a relative entropy
constraint that
bounds the information-processing capacity of the decision-maker.
In the relative entropy we recognize the term $\frac{1}{\beta} \log p(x)$ as our
transformation costs $\rho$ from Theorem~\ref{theo:utility} such that
we can express the negative free energy difference $-\Delta F$ as
\[
-\Delta F = \expect[U] - \expect[R],
\]
where $U(x)=\phi_0(x)-\phi(x)$ represents gains (or losses) and $R(x)
=  \rho(x) - \rho_0(x) $ represents the extra resource costs required
to achieve the gain (or loss) $U$. We can therefore see how the
variational principle of Corollary~\ref{theo:var-principle} formalizes
a trade-off between expected gains (or losses) and information
processing costs.

\section{Summary of Main Concepts}

In decision theory, choices between alternatives are usually formalized as choices between lotteries, where a lottery is formalized as a set $\fs{X}$ of possible outcomes, a probability distribution $p_0$ over $\fs{X}$, and a real-valued function $U$ over $\fs{X}$ called the utility function. In particular expected utility theory predicts that a decision-maker always chooses the lottery with the higher expected utility $\expect[U] = \sum_x p_0(x) U(x)$.
Here we introduce the notion of a \emph{bounded lottery} as a lottery that is additionally characterized by a \emph{resource parameter} $\beta \in \reals$ that captures the resource constraints of the decision-maker.

We have derived a thermodynamic framework for bounded lotteries
from simple axioms that measure information processing cost---see also \cite{Ortega2011}.
The most important difference of bounded decision-making compared to perfectly rational decision-making is that the bounded decision-maker will not be able to choose infallibly
the best lottery. In fact, the resource constraints lead to stochastic choice behavior which
can be characterized by a probability distribution. The decision process then transforms an initial choice probability $p_0$ into a final choice probability $p$ by taking into account the utility gains (or losses) \emph{and} the transformation costs. This transformation process can be formalized as
\begin{equation}
    \label{eq:choice}
    p(x) = \frac{1}{Z} p_0(x) e^{\beta U(x)},
    \quad\text{where}\quad
    Z = \sum_{x'} p_0(x') e^{\beta U(x')}.
\end{equation}
Accordingly, the choice pattern of the decision-maker is predicted by the probability $p$. Crucially, the probability $p$ extremizes the variational principle
\begin{equation}
\label{eq:free-energy}
    \max_{p} \bigg\{
        \sum_x p(x) U(x) - \frac{1}{\beta} \sum_x p(x) \log \frac{p(x)}{p_0(x)}
    \bigg\}.
\end{equation}
These two terms can be interpreted as determinants of bounded rational decision-making in that they formalize a trade-off between an expected utility gain (first term) and the information processing cost of transforming $p_0$ into $p$ (second term). The certainty equivalent value of a bounded lottery can be obtained by inserting the choice probability $p$ from \eqref{eq:choice} into \eqref{eq:free-energy}, yielding
\begin{equation}
    V = \frac{1}{\beta} \log \biggl( \sum_x p_0(x) e^{\beta U(x)} \biggr),
\end{equation}
which corresponds to the log partition sum. For different values of $\beta$, the certainty equivalent takes the following limits
\begin{eqnarray}
\lim_{\beta \rightarrow \infty} \frac{1}{\beta} \log Z &=& \max_x U(x) \nonumber \\
\lim_{\beta \rightarrow 0} \frac{1}{\beta} \log Z &=& \sum_x p_0(x) U(x) \nonumber \\
\lim_{\beta \rightarrow -\infty} \frac{1}{\beta} \log Z &=& \min_x U(x) \nonumber .
\end{eqnarray}
The case $\beta \rightarrow \infty$ corresponds to the perfectly rational actor (or the extremely optimistic observer), the case $\beta \rightarrow -\infty$ corresponds to the perfectly ``anti-rational'' actor (or the extremely pessimistic observer) and the case $\beta \rightarrow 0$ corresponds to the actor that has no resources (or the risk-neutral observer) such that the best one can expect is the expected gain or loss. For illustration see Figure~2.

\begin{figure}[tb]
\label{fig:ceversusbeta}
\begin{center}
    \vspace*{.05in}
    \footnotesize
    \psfrag{l1}[l]{\textbf{a)}}
    \psfrag{l2}[l]{\textbf{b)}}
    \psfrag{a1}[l]{$U_{\max}$}
    \psfrag{a2}[l]{$U_{\min}$}
    \psfrag{a3}[l]{$\expect[U]$}
    \psfrag{a4}[c]{$\beta$}
    \psfrag{b1}[c]{$\beta_1$}
    \psfrag{b2}[c]{$\beta_2$}
    \psfrag{b3}[c]{$\beta_3$}
    \includegraphics[width=8cm]{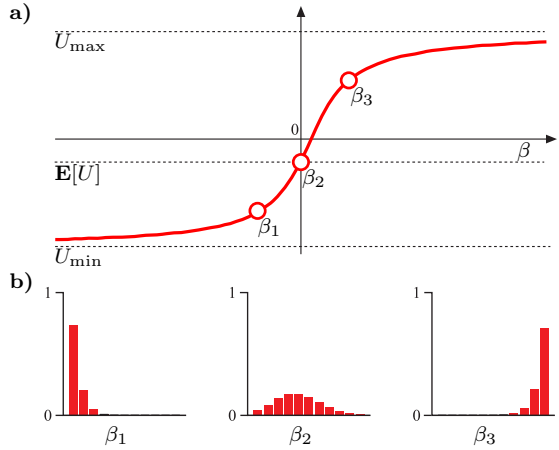}
    \caption{a) Negative free energy difference $\Delta F$ versus the resource parameter~$\beta$. The resource parameter allows modeling decision-makers with bounded resources, either when generating their own actions ($\beta > 0$) or when anticipating their environment ($\beta < 0$). The negative free energy difference corresponds to the certainty equivalent. b) Distribution over the outcomes depending on the resource parameter $\beta$. For large positive $\beta$ the distribution concentrates on the outcome with maximum gain $U_{\max}$. For large negative $\beta$ the distribution concentrates on the worst outcome with gain $U_{\min}$. For $\beta = 0$ the outcomes follow the given distribution $p_0$.}
\end{center}
\end{figure}

\section{Bounded Rationality and Satisficing}

Herbert Simon \cite{Simon1956} proposed in the 50s that bounded rational decision-makers do not commit to an unlimited optimization by searching for the absolute best option. Rather, they follow a strategy of \emph{satisficing}, i.e. they settle for an option that is \emph{good enough} in some sense. Since then, it has been debated whether \emph{satisficing} decision-makers can be described as bounded rational decision-makers that act optimally under resource constraints or whether optimization is the wrong concept altogether \cite{Gigerenzer2001}.
If decision-makers did indeed explicitly attempt to solve such a constrained optimization problem, this would lead to an infinite regress and the paradoxical situation that a bounded rational decision-maker would have to solve a more complex (i.e. constrained) optimization problem than a perfectly rational decision-maker.

To resolve this paradox, the bounded rational decision maker must not be able to reason about his constraints. He just searches randomly for the best option, until his resources run out. An observer will then be able to assign a probability distribution to the decision-maker's choices and investigate how this probability distribution changes depending on the available resources. Consider, for example, an anytime algorithm that will compute a solution more and more precisely the more time it has at its disposal. As one does not want to wait forever for an answer, the anytime computation will be interrupted at some point where one assumes that the answer is going to be good enough. This concept of satisficing can be used to interpret Equation~\ref{eq:probdist} which describes the choice rule of a bounded rational decision-maker.

Consider the problem of picking the largest number in a sequence $U_0, U_1, U_2, \ldots$ of i.i.d. data, where each $U_i \in \fs{U}$ is drawn from a source with probability distribution $\mu$. This could be, for instance, an urn with numbered balls that we draw with replacement and we always keep track of the largest number seen so far. After $m$ draws the largest number will be given by
\[
    v \define \max \{ U_1, U_2, \ldots, U_m \}.
\]
Naturally, the larger the number of draws, the higher the chances of observing a large number. The cumulative distribution function of choosing $v$ after $m$ draws is given by
\begin{equation}
\label{eq:maxdist}
F_{m} (v) = F_0(v)^{m},
\end{equation}
where $F_0$ is the cumulative distribution function of $\mu$ \cite{Gumbel1958}.
If we only cared about finding the maximum with absolute certainty then we would need to draw an infinite amount of times.
However, a bounded rational decision-maker would stop after a certain time, when he feels that the benefit of further exploration does not justify the effort of further drawings. Thus, the number of draws in this example can be regarded as a resource and the numbers on the balls can be regarded as utilities. The behavior of the bounded rational decision-maker is then stochastic even though he acts perfectly deterministically, in the sense that he chooses option $v$ with probability \eqref{eq:maxdist} given the resource constraint $m$.
According to \eqref{eq:maxdist}, the more resources a decision-maker spends, the more he resembles a perfectly rational decision-maker that chooses the maximum number (Figure~1a), since the expected utility increases monotonically with the amount of resources spent (Figure~1b). Importantly, however, note that the marginal increase in the expected utility diminishes with larger effort---hence larger and larger effort pays out less and less in the end. Below we formalize this trade-off.

\begin{figure}[tb]
\begin{center}
    \vspace*{.05in}
    \footnotesize
    \psfrag{l1}[c]{\textbf{a)}}
    \psfrag{l2}[c]{\textbf{b)}}
    \psfrag{b1}[c]{$M=0$}
    \psfrag{b2}[c]{$M=8$}
    \psfrag{b3}[c]{$M=32$}
    \psfrag{b4}[c]{$M=128$}
    \psfrag{a1}[c]{$U_{\max}$}
    \psfrag{a3}[c]{$0$}
    \psfrag{a4}[c]{$M$}
    \psfrag{a5}[c]{$\expect[v]$}
    \psfrag{a6}[c]{$\expect[v] - M \cdot c$}
    \includegraphics[width=8cm]{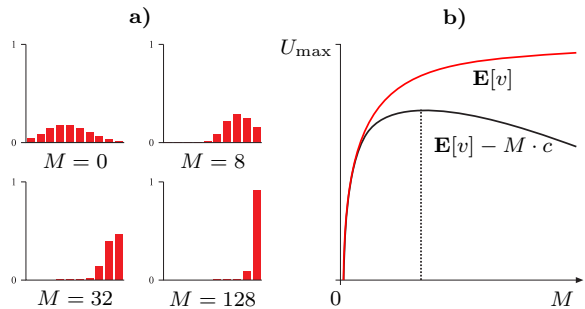}
    \caption{a) Distributions over the maximum for various sample sizes $(M+1)$. The distribution $\mu$ over the ten values $v$ in $\fs{U} = \{1,2,3,\ldots,10\}$ follows a truncated Poisson distribution with parameter $\lambda=5$, as can be seen in the plot for $M=0$. The distribution approaches a delta function over $v = 10$ for increasing values of~$M$. b) The expected maximum $v$ versus sample size $(M+1)$. The incremental gain of the expected maximum is marginally decreasing as the sample size increases (red). If the sampling process is associated with a cost---e.g. $c = 0.02$ per sample in the figure---, then the penalized expected maximum (black) reaches a unique maximum for a finite sample size---the optimal sample size is $M=35$ in the figure.}
    \label{fig:maximum}
\end{center}
\end{figure}

Here we show that the boundedness parameter $\beta$ plays an analogous role to the number of draws $m$. In the limit of a continuous cumulative function $F_0$, the density after $m$ draws is given by $p_m(v) = \frac{d}{dv} F_0(v)^m$. We can now compute the log odds for two random outcomes $v$ and $v'$, which results in
\[
\log \frac{p_m(v)}{p_m(v')} = (m - 1) \log \frac{F_0(v)}{F_0(v')} + \log \frac{\mu(v)}{\mu(v')},
\]
where $F_0(v)$ is again the cumulative of $\mu$. If we require the probabilities $p_m(v)$ to be representable by a distribution of the exponential family such that $p_m(v)=\frac{\mu(v) \exp\left( \alpha U(v) \right)}{\int\,dv' \mu(v') \exp\left( \alpha U(v') \right)}$, we see that the log odds have the following relation
\[
\log \frac{p_m(v)}{p_m(v')} = \alpha \left(U(v) - U(v') \right) + \log \frac{\mu(v)}{\mu(v')}.
\]
We see that $\alpha$ and $m$ play the role of the number of samples or computations.
In general, the following theorem can be shown to hold.

\begin{theorem}\label{theo:temperature-interpretation}
Let $\fs{X}$ be a finite set. Let $Q$ and $M$ be strictly positive probability distributions over $\fs{X}$. Let $\alpha$ be a positive integer. Define $M_\alpha$ as the probability distribution over the maximum of $\alpha$ samples from $M$. Then, there are strictly positive constants $\delta$ and $\xi$ depending only on $M$ such that for all $\alpha$,
\[
    \left|
        \frac{ Q(x) e^{\alpha U(x)} }
        { \sum_{x'} Q(x') e^{\alpha U(x')} }
        - M_\alpha(x)
    \right|
    \leq e^{-(\alpha-\xi) \delta}.
\]
\end{theorem}

Consequently, one can interpret the inverse temperature as a resource parameter that determines how many samples are drawn to estimate the maximum. Note that the distribution~$M$ is arbitrary as long as it has the same support as $Q$. This interpretation can be extended to a negative $\alpha$, by noting that $\alpha U(x) = (-\alpha)(-U(x))$, i.e.\ instead of the maximum we take the minimum of $-\alpha$ samples.

\section{Sequential Decision-Making}

In the case of sequential decision-making the assumption of uniform temperatures has to be relaxed---the proofs of the following theorems can be found in \cite{Ortega2012}. In general, we can then dedicate different amounts of computational resources to each node of a decision tree. However, this requires a translation between a tree with a single temperature and to a tree with different temperatures. This translation can be achieved using the following theorem

\begin{theorem}\label{theo:temp-change}
Let $P$ be the equilibrium distribution for a given inverse temperature~$\alpha$, utility function~$U$ and reference distribution~$Q$. If the temperature changes to~$\beta$ while keeping~$P$ and~$Q$ fixed, then the utility function changes to
\[
    V(x) = U(x) - \Bigl( \tfrac{1}{\alpha} - \tfrac{1}{\beta} \Bigr)
        \log \frac{P(x)}{Q(x)}.
\]
\end{theorem}

If we now define the reward as the change in utility of two subsequent nodes, then the rewards of the resulting decision tree are given by
\begin{eqnarray}\label{eq:rewards}
    R(x_t|x_{<t}) &:=&  \bigl[ V(x_{\leq t}) - V(x_{<t}) \bigr] \nonumber \\
    &=& \bigl[ U(x_{\leq t}) - U(x_{<t}) \bigr]
        - \Bigl( \tfrac{1}{\alpha} - \tfrac{1}{\beta(x_{<t})} \Bigr)
            \log\frac{ P(x_t|x_{<t}) }{ Q(x_t|x_{<t})} \nonumber.
\end{eqnarray}
This allows introducing a collection of node-specific (not necessarily time-specific) inverse temperatures $\beta(x_{<t})$, allowing for a greater degree of flexibility in the representation of information costs. The next theorem states the connection between the free energy and the general decision tree formulation.

\begin{theorem}\label{theo:free-energy-trajectory}
The free energy of the whole trajectory can be rewritten in terms of rewards:
\begin{align}
    \nonumber
    F_\alpha[P] &=
    \sum_{x_{\leq T}} P(x_{\leq T})
    \biggl\{ U(x_{\leq T})
        - \frac{1}{\alpha} \log \frac{P(x_{\leq T})}{Q(x_{\leq T})} \biggr\}
    \\ \label{eq:seq-free-energy}
    &= U(\varepsilon) +
    \sum_{x_{\leq T}} P(x_{\leq T})
        \sum_{t=1}^T \biggl\{
        R(x_t|x_{<t})
        - \frac{1}{\beta(x_{<t})} \log \frac{P(x_t|x_{<t})}{Q(x_t|x_{<t})} \biggr\}.
\end{align}
\end{theorem}

This translation allows applying the free energy principle to each node with a different resource parameter $\beta(x_{<t})$. By writing out the sum in \eqref{eq:seq-free-energy}, one realizes that this free energy has a nested structure where the latest time step forms the innermost variational problem and all other variational problems of the previous time steps can be solved recursively by working backwards in time. This then leads to the following solution:

\begin{theorem}\label{theo:free-energy-solution}
The solution to the free energy in terms of rewards is given by
\[
    P(x_t|x_{<t}) = \frac{1}{Z(x_{<t})}
        Q(x_t|x_{<t}) \exp\Bigl\{
        \beta(x_{<t}) \bigl[ R(x_t|x_{<t})
            + \frac{1}{\beta(x_{\leq t}) } \log Z(x_{\leq t})
            \bigr] \Bigr\},
\]
where $Z(x_{\leq T}) = 1$ and where for all $t<T$
\[
    Z(x_{<t}) = \sum_{x_t} Q(x_t|x_{<t}) \exp\Bigl\{
        \beta(x_{<t}) \bigl[ R(x_t|x_{<t})
            + \frac{1}{\beta(x_{\leq t}) } \log Z(x_{\leq t})
            \bigr] \Bigr\}.
\]
\end{theorem}

\section{Limit Cases of Bounded Rational Control}

\begin{figure}[tb]
\begin{center}
    \vspace*{.05in}
    \footnotesize
    \psfrag{l1}[c]{$\beta_{1}$}
    \psfrag{l2}[c]{$\beta_{2}$}
    \psfrag{v1}[c]{$+\infty$}
    \psfrag{v2}[c]{$-\infty$}
    \psfrag{v3}[c]{$+\infty$}
    \psfrag{v4}[c]{$-\infty$}
    \psfrag{a1}[c]{Anti-}
    \psfrag{a2}[c]{Rational}
    \psfrag{a3}[c]{Irrational}
    \psfrag{a4}[c]{Rational}
    \psfrag{b1}[l]{Risk-}
    \psfrag{b2}[l]{Seeking}
    \psfrag{b3}[l]{Risk-}
    \psfrag{b4}[l]{Neutral}
    \psfrag{b5}[l]{Risk-}
    \psfrag{b6}[l]{Averse}
    \includegraphics[width=7cm]{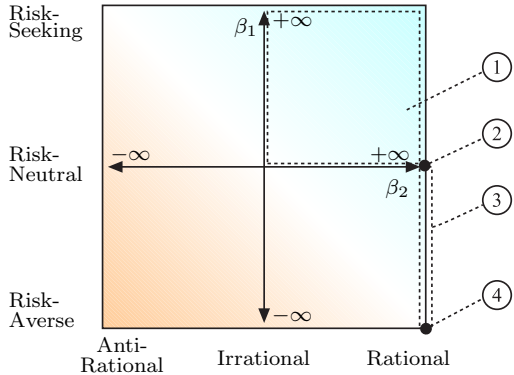}
    \caption{Schematic illustration of how resource parameters can model a range of decision-making schemes: (1)- risk-seeking, bounded rational; (2) risk-neutral, perfectly rational; (3) risk-averse, perfectly rational; and (4) robust, perfectly rational.}
    \label{fig:betas}
\end{center}
\end{figure}

As described in the previous section, the belief and action probabilities of an agent in a sequential decision-making setup can be determined by recursion of the log-partition function
\begin{equation}
\label{eq:equation2}
    V(x_{<t})
    = \frac{1}{\beta(x_{<t})} \log\biggl\{
        \sum_{x_t} Q(x_t|x_{<t}) \exp\Bigl\{
        \beta(x_{<t}) \bigl[ R(x_t|x_{<t})
            + V(x_{\leq t})
            \bigr] \Bigr\} \biggr\},
\end{equation}
where we have introduced $V(x_{\leq t})=\frac{1}{\beta(x_{\leq t}) } \log Z(x_{\leq t})$.
If $x_t$ is an action variable then $Q(x_t|x_{<t})$ reflects the prior policy and the agent's rationality $\beta(x_{<t})$ determines in how far the value $R(x_t|x_{<t})+V(x_{\leq t})$ can be optimized by the agent. If $x_t$ is an observation variable then $Q(x_t|x_{<t})$ reflects the agent's prior belief and the rationality of the environment $\beta(x_{<t})$ indicates how much one should deviate from the prior belief considering the possible values $R(x_t|x_{<t})+V(x_{\leq t})$. Depending on $\beta(x_{<t})$, different decision-making schemes can be recovered---compare Figure~3.

\begin{enumerate}

\item{\textbf{KL control.}}
When assuming a history-independent loss function $r(x_t)$, Markov probabilities $p_0(x_t|x_{t-1})$ and $\beta(x_{<t})=\beta$ for all $x_{<t}$, Equation~\eqref{eq:equation2} simplifies to a recursion that is equivalent to $z$-iteration which has previously been suggested in \cite{Todorov2006,Todorov2009} to approximately solve MDPs by means of linear algebra---see~\cite{BraunOrtega2011,OrtegaBraun2011b} for details of this equivalence relation. In \cite{Todorov2006,Todorov2009}
the transition probabilities of the MDP are controlled directly and the control costs are given by the Kullback-Leiber divergence of the manipulated state transition
probabilities with respect to a baseline distribution that describes the passive dynamics. In our framework, this kind of KL control corresponds to the special case where all random variables are action variables and the agent has boundedness parameter $\beta$. The stochasticity in this case, however, is not thought to arise from environmental passive dynamics, but rather is a direct consequence of bounded rational control in a (possibly) deterministic environment. The continuous case of KL control relies on the formalism of path integrals \cite{Kappen2005,Theodorou2010}, but essentially the same relation to bounded rationality can be established---see~\cite{BraunOrtega2011} for details.

\item{\textbf{Optimal stochastic control.}}
When assuming $\beta(x_{<t})\rightarrow\infty$ for all action variables and $\beta(x_{<t})\rightarrow 0$ for all observation variables, we approach the limit of the
perfectly rational decision-maker in a stochastic environment. In this limit, the log-partition function converges to the expected utility and the decision-maker acts deterministically so as to maximize the expected utility. For action variables, recursion~\eqref{eq:equation2} becomes the well-known
\emph{Bellman Optimality Equation} \cite{Bellman1957}---see \cite{BraunOrtega2011,OrtegaBraun2011b} for details.

\item{\textbf{Risk-sensitive control.}}
Risk-sensitive control \cite{Whittle1990}
corresponds to a decision-maker with $\beta(x_{<t})\rightarrow\infty$ for all action variables and $\beta(x_{<t}) \neq 0$ for observation variables. Risk-sensitivity in the context of continuous KL control has been previously proposed in \cite{Broek2010}. Mean-variance decision criteria used in finance can be equally derived \cite{Markowitz1952}. Risk-sensitive decision-makers do not simply maximize the expectation of the utility, but also consider higher-order cumulants by optimizing a \emph{stress function} given by the log partition sum.
A risk-averse decision-maker ($\beta(x_{<t}) < 0$), for example, discounts variability off the expected utility.
In contrast, risk-seeking decision-makers ($\beta(x_{<t}) > 0$) add value to the expected utility in the face of variability. Risk-sensitivity biases the beliefs about the environment optimistically (collaborative environment) or pessimistically (adversarial environment). Alternatively, one could regard a collaborative environment also as a bounded rational controller that can choose its own observation---that is the environment behaves like an extension of the agent with partial control.
Importantly, the stress function is typically assumed in risk-sensitive control
schemes in the literature, whereas here it falls out naturally---see~\cite{OrtegaBraun2011b} for more details.

\item{\textbf{Robust control.}}
When assuming $\beta(x_{<t})\rightarrow\infty$ for all action variables and $\beta(x_{<t})\rightarrow -\infty$ for all observation variables, we approach the limit of the
robust decision-maker in an unknown environment.
When $\beta(x_{<t}) \rightarrow -\infty$, the decision-maker makes a worst case assumption about the environment, namely that it is strictly adversarial and perfectly rational. This leads to the well-known game-theoretic minimax problem. Minimax problems have been used to reformulate robust control problems that allow controllers to cope with model uncertainties \cite{Basar1991,Hansen2008}. Robust control problems have long been known to be related
to risk-sensitive control \cite{Jacobson1973,Boyle1988}. Here we derived both control types from the same variational principle---see~\cite{OrtegaBraun2011b} for more details.

\end{enumerate}

\section{Discussion}

In the proposed thermodynamic interpretation of bounded rationality,
agents with limited resources search for a maximum over a set by randomly drawing elements
from this set. This random search leads to a utility function that is marginally decreasing
when more search effort is allocated. When such agents pay a search cost, the bounded
rational optimum is to abort the search as soon as the marginal returns are equal to the
search cost. The resulting trade-off between utility maximization and resource costs can be quantified by the Kullback-Leibler divergence with respect to an initial policy or belief.
This initial probability distribution corresponds to the initial state of a thermodynamic system
that changes when a new potential is imposed. The difference in the potential
corresponds to utility gains or losses in economic choice. The difference in the free energy corresponds to physical work and the economic certainty-equivalent. Thus, gains or losses that
are associated with uncertainty are effectively devalued or overvalued, depending on the sign
of the resource parameter. This way risk-sensitivity, robustness to model uncertainty and
game-theoretic minimax-strategies can arise naturally.

\paragraph{Bounded rationality.}
Starting with Simon \cite{Simon1956}, bounded rationality has been extensively studied
in psychology, economics, political science, industrial organization, computer science and artificial intelligence research---see for example \cite{Lipman1995,Russell1995a,Russell1995b,Aumann1997,Rubinstein1998,Gigerenzer2001,Kahnemann2003,Spiegler2011}.
Additionally, numerous experiments in behavioral economics have shown that humans systematically violate perfect rationality, that is they are bounded rational \cite{Camerer2003}.
Probably the most closely related approach to bounded rationality with respect to the present article is quantal response equilibrium (QRE) game theory \cite{QRE1995,QRE1998,Anderson2002,Wolpert2012}.
QRE models assume bounded rational players whose choice probabilities are given by the Boltzmann distribution and whose rationality is determined by a temperature parameter. Interactions of such bounded rational players can lead to game-theoretic solutions that deviate from the Nash equilibrium. The QRE model is a special case of the model presented here where all prior probabilities are assumed to be uniform. These prior probabilities are crucial when defining the certainty-equivalent that ranges from minimum to maximum via the expected utility. As the certainty-equivalent corresponds to physical work, this also allows to relate bounded rational decision-making to thermodynamic processes. The distinction of a prior policy and a utility that is optimized to some extent is  fundamental to the notion of bounded rationality proposed in this paper and therefore also affords a qualitative advance of the bounded rationality model in QRE models.

\paragraph{Information theory in control and game theory.}
As already discussed, a number of papers have suggested the use of the relative entropy as a cost function for control \cite{Todorov2006,Todorov2009,Peters2010,Kappen2012}. Previously, Saridis \cite{Saridis1988} has framed optimal and adaptive control as entropy minimization problems.
Statistical physics has also served as an inspiration to a number of other studies, for example, to an information-theoretic
approach to interactive learning \cite{Still2009}, to use information theory to approximate joint strategies in games with bounded rational
players \cite{Wolpert2004} and to the problem of optimal search \cite{Stone1975,Jaynes1985}, where the utility losses correspond directly to search effort. Recently, Tishby and Polani \cite{Tishby2011} have shown how to apply information theory to understand information flow in the action-observation cycle.
The contribution of our study is to devise information-theoretic axioms
to quantify search costs in bounded optimization problems.
This allows for a unified treatment of control and game-theoretic problems, as well as estimation and learning problems for both perfectly rational and bounded rational agents.
In the future it will be interesting to relate the thermodynamic resource costs of
bounded rational agents to more traditional notions of resource
costs in computer science like space and time requirements of algorithms \cite{Vitanyi2005}.

\paragraph{Variational Preferences.}
In the economic literature the Kullback-Leibler divergence has appeared in the context of
multiplier preference models that can deal with model uncertainty \cite{Hansen2008}.
Especially, it has been proposed that a bound on the Kullback-Leibler divergence could be used to indicate how much of a deviation from a proposed model $p_0$ is allowed when computing robust decision strategies that work under a range of models in the neighborhood of $p_0$. In variational preference models \cite{Rustichini2006} this is generalized to models of the form
\[
f \succeq g \Longleftrightarrow \min_p \left( \int u(f) dp + c(p) \right) \geq \min_p \left( \int u(g) dp + c(p) \right),
\]
where $c(p)$ can be interpreted as an ambiguity index that can explain effects
of ambiguity aversion. The thermodynamic certainty-equivalent of work---computed as the log-partition sum---also falls within this preference model. However, an important difference is that the
choice in a thermodynamic system is not deterministic with respect to the certainty-equivalent,
but stochastic following a generalized Boltzmann distribution. Due to this stochasticity of the
choice behavior itself, the thermodynamic model can be linked to both bounded rationality and
model uncertainty, whereas variational preference models have so far concentrated on explaining
effects of ambiguity aversion and model uncertainty.

\paragraph{Ellsberg's and Allais' paradox.}
Two of the most famous deviations from expected utility theory that have been consistently observed in human decision-making are the paradoxa of Ellsberg \cite{Ellsberg1961} and Allais \cite{Allais1953}. While the first paradox has encouraged a large literature dealing with model uncertainty \cite{Hansen2008}, the latter paradox has led to the development of prospect theory \cite{Kahneman1979,Tversky1992}.
Ellsberg could show that human choice in the face of ambiguity differs from decision-making under risk where precise probability models are available. Humans typically tend to avoid ambiguous options, rather than choosing the option with higher expected utility. The observed ambiguity aversion can be modeled straightforwardly by a bounded rational decision-maker by allowing some degree of minimaxing in the spirit of a risk-sensitive controller---see Supplementary Material for details. Allais could show that humans frequently reverse their preferences in choice tasks that may not lead to preference reversals according to expected utility theory. These reversals typically occur for different levels of riskiness of the same choices. The explanation of the Allais paradox within the framework of bounded rationality is not as straightforward as the Ellsberg paradox, but may involve context-dependent changes of the boundedness parameter or biases in the decision-making process that lead to a \emph{generalized quasi-linear mean model} \cite{Nagumo1930,Kolmogorov1930,Finetti1931,Hong1983}, which provides an alternative account of preference reversals of the Allais type without violating the principle of stochastic dominance---see Supplementary Material for more details.

\paragraph{Stochastic Choice.}
Stochastic choice rules have been extensively studied in the psychological and econometric literature, in particular logit choice models based on the Boltzmann distribution \cite{Luce2000,Train2009}. The literature on Boltzmann distributions for decision-making goes back to Luce \cite{Luce1959}, extending through McFadden \cite{McFadden1974}, Meginnis \cite{Meginnis1976}, Fudenberg \cite{Fudenberg1993} and Wolpert \cite{Wolpert2004,Wolpert2011,Wolpert2012}.
Luce \cite{Luce1959} has studied stochastic choice rules of the form $p(x_i) \sim \frac{w_i}{\sum_j w_j}$, which includes the Boltzmann distribution and the ``softmax''-rule known in the reinforcement learning literature \cite{Sutton1998}.
McFadden \cite{McFadden1974} has shown that such distributions can arise, for example, when utilities are contaminated
with additive noise following an extreme value distribution. While stochastic choice models
are generally accepted to  account for human choices better than their deterministic counterparts \cite{Rieskamp2008,gläscher2010sta,mcdannald2012mod}, they have also been strongly criticized, especially for a property known as \emph{independence of irrelevant alternatives} (IIA). Similar to the independence axiom in expected utility theory, IIA implies that the ratio of two choice probabilities does not depend on the presence of a third irrelevant alternative in the choice set. What distinguishes the free energy equations from above choice rules is that stochastic choice behavior is described by a generalized exponential family distribution of the form $p(x) \sim p_0(x) \exp(\beta U(x))$. Changing the choice set might in general also change the prior $p_0(x)$,
but more importantly it might also change the resource parameter $\beta$.

\paragraph{Diffusion-to-bound models.}
Diffusion-to-bound models typically model the process of binary decision-making as a random walk process that terminates once it hits one of two given decision bounds \cite{BusemeyerD02}. Each time step of the random walk provides noisy evidence towards one of the two options. This implements a natural speed-accuracy trade-off: the further away the bounds the more reliable the decision will be, as the noise can be averaged out, but also the longer one has to wait. The resulting choice probabilities are identical to the choice probabilities of a bounded rational decision-maker
if we relate the decision bound of the random walk with the boundedness parameter in \ref{eq:probdist}---see Supplementary Material for details. The boundedness parameter can then also be shown to be proportional to the time required for the decision-making process.
Decision-to-bound models have been widely used in behavioral psychology and neuroscience
to explain probabilistic choice and reaction times in psychometric experiments---see \cite{Bogacz2006} for a review. Decision-makers that apply the decision-to-bound model may be regarded as bounded rational decision-makers from a normative point of view.

\paragraph{Free Energy Principle.}
A central property of closed thermodynamic systems is that they minimize free energy. A free energy principle based on the variational Bayes approach has recently also been proposed as a theoretical framework to understand brain function \cite{Friston2009,Friston2010}. In this framework generative models of the form $p(y|h,a)$ explain how hidden causes $h$ in the environment and actions $a$ produce observations $y$. The brain uses an approximative distribution $Q(h;a)$ to determine the hidden causes. The free energy
\[
F = - \int dh \, Q(h;a) \ln P(y,h|a) - \int dh \, Q(h;a) \ln Q(h;a)
\]
measures how well the brain is doing with this approximation. According to \cite{Friston2009,Friston2010}, action and perception consist in choosing $a$ and $Q$ respectively so as to minimize this free energy.
In light of the thermodynamic view of free energy, maximizing the likelihood $- \ln P(y,h|a)$---or minimizing surprise---is a particular choice of potential function $\phi$, where the boundedness consists in being restricted to model class $Q$ instead of having full disposal of $p(y|h,a)$. More generally, variational Bayes methods that use particular classes of distributions to approximate the posterior could thus be regarded as a form of bounded inference within this picture.

\section{Conclusion}
Thermodynamics provides a framework for bounded rationality that can be both descriptive and prescriptive. It is descriptive in the sense that it describes behavior that is clearly sub-optimal from the point of view of a perfect rational decision-maker with infinite resources. It is prescriptive in the sense that it prescribes how a bounded rational actor should behave optimally given resource constraints formalized by $\beta$. As we have argued in this paper, bounded rational decision-making provides an overarching principle in both senses in economics, engineering, artificial intelligence, psychology and neuroscience.

A thermodynamic model of bounded rational decision-making has also two advantages over traditional decision theory of perfect rationality. First, it allows connecting computational processes with real physical processes, for example how much entropy they generate and how much energy they require \cite{Tribus1971}. Second, it suggests a notion of intelligence that is closely related to the process of evolution. It is straightforward to see that bounded rational controllers of the form~\eqref{eq:choice} share their structure with Bayes' rule,
where we identify the prior $p_0(x)$, the likelihood model $e^{\beta U(x)}$ and the posterior $p(x)$, normalized by the partition function, thus, establishing a close link between inference and control
\cite{OrtegaBraun2010e}. Furthermore, both bounded rational controllers and Bayes' rule share their structural form with discrete replicator dynamics that model evolutionary processes \cite{Shahlizi2009},
where samples (a population) are pushed through a fitness function (likelihood, gain function) that biases the distribution of the population, thereby transforming a distribution $p_0$ to a new distribution $p$. In this picture different hypotheses $x$ compete for probability mass over subsequent iterations, favoring those $x$ that have a lower-than-average cost.
Just like the evolutionary random processes of variation and selection created intelligent organisms on a phylogenetic timescale, similar random processes might underlie (bounded) intelligent behavior in individuals on an ontogenetic timescale.

%%%%%%%%%%%%%%%%%%%%%%%%%%%%%%%%%%%%%%%%%%%%%%%%%%%%%%%%%%%%%%%%%%%%
% Acknowledgements                                                 %
%%%%%%%%%%%%%%%%%%%%%%%%%%%%%%%%%%%%%%%%%%%%%%%%%%%%%%%%%%%%%%%%%%%%

\section{Acknowledgments}
This study was supported by the DFG, Emmy Noether grant BR4164/1-1.

%%%%%%%%%%%%%%%%%%%%%%%%%%%%%%%%%%%%%%%%%%%%%%%%%%%%%%%%%%%%%%%%%%%%
% Bibliography                                                     %
%%%%%%%%%%%%%%%%%%%%%%%%%%%%%%%%%%%%%%%%%%%%%%%%%%%%%%%%%%%%%%%%%%%%

%GATHER{bibliography.bib}
\bibliographystyle{vancouver}
\bibliography{bibliography}

\begin{thebibliography}{10}

\bibitem{Gigerenzer1999}
Gigerenzer G, Todd PM, {ABC Research Group}.
\newblock {Simple Heuristics That Make Us Smart}.
\newblock New York: Oxford University Press; 1999.

\bibitem{Gladwell2005}
Gladwell M.
\newblock Blink: the power of thinking without thinking.
\newblock New York: Little, Brown and Company; 2005.

\bibitem{NivDawJoelEtAl07}
Niv Y, Daw ND, Joel D, Dayan P.
\newblock Tonic dopamine: opportunity costs and the control of response vigor.
\newblock Psychopharmacology (Berl). 2007 Apr;191(3):507--520.

\bibitem{daw2012}
Daw ND.
\newblock `Model-based reinforcement learning as cognitive search:
  neurocomputational theories'.
\newblock In: Cognitive search: evolution algorithms and the brain. Boston: MIT
  Press; 2012. .

\bibitem{Neumann1944}
Von~Neumann J, Morgenstern O.
\newblock Theory of Games and Economic Behavior.
\newblock Princeton: Princeton University Press; 1944.

\bibitem{Savage1954}
Savage LJ.
\newblock The Foundations of Statistics.
\newblock New York: John Wiley and Sons; 1954.

\bibitem{Fishburn1982}
Fishburn P.
\newblock The Foundations of Expected Utility.
\newblock Dordrecht: D. Reidel Publishing; 1982.

\bibitem{Simon1972}
Simon H.
\newblock Theories of Bounded Rationality.
\newblock In: Radner CB, Radner R, editors. Decision and Organization.
  Amsterdam: North Holland Publ.; 1972. p. 161--176.

\bibitem{Simon1984}
Simon H.
\newblock Models of Bounded Rationality.
\newblock Cambridge, MA: {MIT} Press; 1984.

\bibitem{Rubinstein1998}
Rubinstein A.
\newblock Modeling Bounded Rationality.
\newblock Cambridge, MA: {MIT} Press; 1999.

\bibitem{Gigerenzer2001}
Gigerenzer G, Selten R.
\newblock In: Bounded rationality: the adaptive toolbox. Cambridge, MA: {MIT}
  Press; 2001. .

\bibitem{Feynman1996}
Feynman RP.
\newblock The Feynman Lectures on Computation.
\newblock Addison-Wesley; 1996.

\bibitem{Tribus1971}
Tribus M, McIrvine EC.
\newblock Energy and Information.
\newblock Scientific American. 1971;225:179--188.

\bibitem{Landauer1961}
Landauer R.
\newblock Irreversibility and Heat Generation in the Computing Process.
\newblock IBM Journal of Research and Development. 1961;5(3):183--191.

\bibitem{Bennett1973}
Bennett CH.
\newblock Logical Reversibility of Computation.
\newblock IBM Journal of Research and Development. 1973;17(6):525--532.

\bibitem{Bennett1982}
Bennett CH.
\newblock The thermodynamics of computation—a review.
\newblock International Journal of Theoretical Physics. 1982;21(12):905--–940.

\bibitem{Mackay2003}
MacKay DJC.
\newblock Information Theory, Inference, and Learning Algorithms.
\newblock Cambridge University Press; 2003.

\bibitem{OrtegaBraun2010b}
Ortega PA, Braun DA.
\newblock A conversion between utility and information.
\newblock In: Proceedings of the third conference on artificial general
  intelligence. Atlantis Press; 2010. p. 115--120.

\bibitem{Ortega2011}
Ortega PA.
\newblock A Unified Framework for Resource-Bounded Autonomous Agents
  Interacting with Unknown Environments.
\newblock Department of Engineering, University of Cambridge, UK; 2011.

\bibitem{Shannon1948}
Shannon CE.
\newblock A mathematical theory of communication.
\newblock Bell System Technical Journal. 1948 Jul and Oct;27:379--423 and
  623--656.

\bibitem{Csiszar2008}
Csisz\'{a}r I.
\newblock Axiomatic Characterizations of Information Measures.
\newblock Entropy. 2008;10:261--273.

\bibitem{Callen1985}
Callen HB.
\newblock Thermodynamics and an introduction to thermostatistics.
\newblock New York: John Wiley \& Sons; 1985.

\bibitem{Simon1956}
Simon HA.
\newblock Rational choice and the structure of the environment.
\newblock Psychological Review. 1956;63(2):129--38.

\bibitem{Gumbel1958}
Gumbel EJ.
\newblock Statistics of Extremes.
\newblock New York: Columbia University Press; 1958.

\bibitem{Ortega2012}
Ortega PA, A BD.
\newblock Free Energy and the Generalized Optimality Equations for Sequential
  Decision Making.
\newblock arXiv:12053997v1 (presented at the European Workshop for
  Reinforcement Learning). 2012;.

\bibitem{Todorov2006}
Todorov E.
\newblock Linearly solvable Markov decision problems.
\newblock In: Advances in Neural Information Processing Systems. vol.~19; 2006.
  p. 1369--1376.

\bibitem{Todorov2009}
Todorov E.
\newblock Efficient computation of optimal actions.
\newblock Proceedings of the National Academy of Sciences USA.
  2009;106:11478--11483.

\bibitem{BraunOrtega2011}
Braun DA, Ortega PA.
\newblock Path integral control and bounded rationality.
\newblock In: IEEE Symposium on adaptive dynamic programming and reinforcement
  learning; 2011. p. 202--209.

\bibitem{OrtegaBraun2011b}
Ortega PA, Braun DA.
\newblock Information, utility and bounded rationality.
\newblock In: Lecture notes on artificial intelligence. vol. 6830; 2011. p.
  269--274.

\bibitem{Kappen2005}
Kappen HJ.
\newblock A linear theory for control of non-linear stochastic systems.
\newblock Physical Review Letters. 2005;95:200201.

\bibitem{Theodorou2010}
Theodorou E, Buchli J, Schaal S.
\newblock A generalized path integral approach to reinforcement learning.
\newblock Journal of Machine Learning Research. 2010;11:3137--3181.

\bibitem{Bellman1957}
Bellman RE. Dynamic Programming. Princeton, NJ: Princeton University Press;
  1957.

\bibitem{Whittle1990}
Whittle P.
\newblock Risk-sensitive optimal control.
\newblock New York: John Wiley and Sons; 1990.

\bibitem{Broek2010}
van~den Broek JL, Wiegerinck WAJJ, Kappen HJ.
\newblock Risk-sensitive path integral control.
\newblock In: UAI. vol.~6; 2010. p. 1--8.

\bibitem{Markowitz1952}
Markowitz H.
\newblock Portfolio Selection.
\newblock The Journal of Finance. 1952;7:77--–91.

\bibitem{Basar1991}
Ba\c{s}ar T, Bernhard P.
\newblock H-infinity optimal control and related minimax-design problems: a
  dynamic game approach.
\newblock Boston: Birkh\"{a}user; 1991.

\bibitem{Hansen2008}
Hansen LP, Sargent TJ.
\newblock Robustness.
\newblock Princeton: Princeton University Press; 2008.

\bibitem{Jacobson1973}
Jacobson D.
\newblock Optimal stochastic linear systems with exponential performance
  criteria and their relation to deterministic differential games.
\newblock IEEE T Automat Contr AC. 1973;18:124--131.

\bibitem{Boyle1988}
Glover K, Boyle J.
\newblock State-space formulae for all stabilizing controllers that satisfy an
  H-norm bound and relations to risk sensitivity.
\newblock Syst Control Lett. 1988;11:167--172.

\bibitem{Lipman1995}
Lipman B.
\newblock Information Processing and Bounded Rationality: A Survey.
\newblock Canadian Journal of Economics. 1995;28(1):42--67.

\bibitem{Russell1995a}
Russell SJ.
\newblock Rationality and Intelligence.
\newblock In: Mellish C, editor. Proceedings of the Fourteenth International
  Joint Conference on Artificial Intelligence. San Francisco: Morgan Kaufmann;
  1995. p. 950--957.

\bibitem{Russell1995b}
Russell SJ, Subramanian D.
\newblock Provably bounded-optimal agents.
\newblock Journal of Artificial Intelligence Research. 1995;3:575--609.

\bibitem{Aumann1997}
Aumann RJ.
\newblock Rationality and Bounded Rationality.
\newblock Games and Economic Behavior. 1997 October;21(1-2):2--14.

\bibitem{Kahnemann2003}
Kahneman D.
\newblock Maps of Bounded Rationality: Psychology for Behavioral Economics.
\newblock American Economic Review. 2003 December;93(5):1449--1475.

\bibitem{Spiegler2011}
Spiegler R.
\newblock Bounded Rationality and Industrial Organization.
\newblock Oxford: Oxford University Press; 2011.

\bibitem{Camerer2003}
Camerer C.
\newblock Behavioral Game Theory: Experiments in Strategic Interaction.
\newblock Princeton: Princeton University Press; 2003.

\bibitem{QRE1995}
D MR, R PT.
\newblock Quantal Response Equilibria for Normal Form Games.
\newblock Games and Economic Behavior. 1995 July;10(1):6--38.

\bibitem{QRE1998}
Mckelvey R, Palfrey TR.
\newblock {Quantal Response Equilibria for Extensive Form Games}.
\newblock Experimental Economics. 1998;1:9--41.

\bibitem{Anderson2002}
Anderson SP, Goeree JK, Holt CA.
\newblock The logit equilibrium: a perspective on intuitive behavioral
  anomalies.
\newblock Southern Economic Journal. 2002;69:21--47.

\bibitem{Wolpert2012}
Wolpert DH, Harre M, Bertschinger N, Olbrich E, Jost J.
\newblock Hysteresis effects of changing parameters of noncooperative games.
\newblock Physical Review E. 2012;85:036102.

\bibitem{Peters2010}
Peters J, M\"{u}lling K, Altun Y.
\newblock Relative entropy policy search.
\newblock In: AAAI; 2010. .

\bibitem{Kappen2012}
Kappen HJ, G\'{o}mez V, Opper M.
\newblock Optimal control as a graphical model inference problem.
\newblock Machine Learning. 2012;1:1--11.

\bibitem{Saridis1988}
Saridis G.
\newblock Entropy formulation for optimal and adaptive control.
\newblock IEEE Transactions on Automatic Control. 1988;33:713--721.

\bibitem{Still2009}
Still S.
\newblock An information-theoretic approach to interactive learning.
\newblock Europhysics Letters. 2009;85:28005.

\bibitem{Wolpert2004}
Wolpert DH.
\newblock Information theory - the bridge connecting bounded rational game
  theory and statistical physics.
\newblock In: Braha D, Bar-Yam Y, editors. Complex Engineering Systems. Perseus
  Books; 2004. .

\bibitem{Stone1975}
Stone LD.
\newblock Theory of optimal search.
\newblock New York: Academic Press; 1998.

\bibitem{Jaynes1985}
Jaynes ET.
\newblock Entropy and search theory.
\newblock In: Smith CR, Grandy WT, editors. Maximum entropy and Bayesian
  methods in inverse problems. Dordrecht: Reidel; 1985. .

\bibitem{Tishby2011}
Tishby N, Polani D.
\newblock Information Theory of Decisions and Actions.
\newblock In: Vassilis T Hussain, editor. Perception-reason-action cycle:
  Models, algorithms and systems. Berlin: Springer; 2011. .

\bibitem{Vitanyi2005}
Vitanyi PMB.
\newblock Time, space, and energy in reversible computing.
\newblock In: Proceedings of the 2nd ACM conference on Computing frontiers;
  2005. p. 435--444.

\bibitem{Rustichini2006}
Maccheroni F, Marinacci M, Rustichini A.
\newblock Ambiguity aversion, robustness, and the variational representation of
  preferences.
\newblock Econometrica. 2006;74:1447--1498.

\bibitem{Ellsberg1961}
Ellsberg D.
\newblock Risk, Ambiguity and the Savage Axioms.
\newblock The Quaterly Journal of Economics. 1961;75:643--669.

\bibitem{Allais1953}
Allais M.
\newblock Le comportment de l'homme rationnel devant la risque: critique des
  postulats et axiomes de l'ecole Americaine.
\newblock Econometrica. 1953;21:503--546.

\bibitem{Kahneman1979}
Kahneman D, Tversky A.
\newblock Prospect Theory: An Analysis of Decision under Risk.
\newblock Econometrica. 1979;47:263--291.

\bibitem{Tversky1992}
Tversky A, Kahneman D.
\newblock Advances in prospect theory: Cumulative representation of
  uncertainty".
\newblock Journal of Risk and Uncertainty. 1992;5:297--323.

\bibitem{Nagumo1930}
Nagumo M.
\newblock \"{U}ber eine Klasse der Mittelwerte.
\newblock Japan Journal of Mathematics. 1930;7:71--79.

\bibitem{Kolmogorov1930}
Kolmogorov A.
\newblock Sur la notion de la moyenne.
\newblock Rendiconti accademia dei lincei. 1930;12:388--391.

\bibitem{Finetti1931}
de~Finetti B.
\newblock Sul concetto di media.
\newblock Giornale dell' istituto italiano degli attuari. 1931;2:369--396.

\bibitem{Hong1983}
Hong CS.
\newblock A generalization of the quasilinear mean with application to the
  measurement of income inequality and decision theory resolving the Allais
  paradox.
\newblock Econometrica. 1983;51:1065--1092.

\bibitem{Luce2000}
Luce RD.
\newblock Utility of gains and losses: measurement-theoretical and experimental
  approaches.
\newblock Mahwah, NJ: Erlbaum; 2000.

\bibitem{Train2009}
Train KE.
\newblock Discrete Choice Methods with Simulation.
\newblock 2nd ed. Cambridge: Cambridge University Press; 2009.

\bibitem{Luce1959}
Luce RD.
\newblock Individual choice behavior.
\newblock Oxford: Wiley; 1959.

\bibitem{McFadden1974}
McFadden D.
\newblock Conditional logit analysis of qualitative choice behavior.
\newblock In: Zarembka P, editor. Frontiers in econometrics. New York: Academic
  Press; 1974. .

\bibitem{Meginnis1976}
Meginnis JR.
\newblock A new class of symmetric utility rules for gambles, subjective
  marginal probability functions, and a generalized Bayes' rule.
\newblock Proceedings of the American Statistical Association, Business and
  Economic Statistics Section. 1976;p. 471--476.

\bibitem{Fudenberg1993}
Fudenberg D, Kreps D.
\newblock Learning mixed equilibria.
\newblock Games and Economic Behavior. 1993;5:320--367.

\bibitem{Wolpert2011}
Lee R, Wolpert DH.
\newblock Game-Theoretic Modeling of Human Behavior in Mid-Air Collisions.
\newblock In: T~Guy MK, Wolpert DH, editors. Decision-Making with Imperfect
  Decision Makers. Springer; 2011. .

\bibitem{Sutton1998}
Sutton RS, Barto AG.
\newblock Reinforcement Learning: An Introduction.
\newblock Cambridge, MA: MIT Press; 1998.

\bibitem{Rieskamp2008}
Rieskamp J.
\newblock The probabilistic nature of preferential choice.
\newblock Journal of Experimental Psychology: Learning, Memory, and Cognition.
  2008;34:1446--1465.

\bibitem{gläscher2010sta}
Gläscher J, Daw N, Dayan P, O'Doherty JP.
\newblock States versus rewards: dissociable neural prediction error signals
  underlying model-based and model-free reinforcement learning.
\newblock Neuron. 2010;66(4):585--95.

\bibitem{mcdannald2012mod}
McDannald MA, Takahashi YK, Lopatina N, Pietras BW, Jones JL, Schoenbaum G.
\newblock Model-based learning and the contribution of the orbitofrontal cortex
  to the model-free world.
\newblock Eur J Neurosci. 2012;35(7):991--6.

\bibitem{BusemeyerD02}
Busemeyer JR, Diederich A.
\newblock Survey of decision field theory.
\newblock Mathematical Social Sciences. 2002;43(3):345--370.

\bibitem{Bogacz2006}
Bogacz R, Brown E, Moehlis J, Holmes P, Cohen JD.
\newblock The physics of optimal decision making: a formal analysis of models
  of performance in two-alternative forced-choice tasks.
\newblock Psychological Review. 2006;113:700--765.

\bibitem{Friston2009}
Friston K.
\newblock The free-energy principle: a rough guide to the brain?
\newblock Trends in Cognitive Science. 2009;13:293--301.

\bibitem{Friston2010}
Friston K.
\newblock The free-energy principle: a unified brain theory?
\newblock Nature Review Neuroscience. 2010;11:127--138.

\bibitem{OrtegaBraun2010e}
Ortega PA, Braun DA.
\newblock A minimum relative entropy principle for learning and acting.
\newblock Journal of Artificial Intelligence Research. 2010;38:475--511.

\bibitem{Shahlizi2009}
Shahlizi CR.
\newblock Dynamics of bayesian updating with dependent data and misspecified
  models.
\newblock Electronic Journal of Statistics. 2009;3:1039–--1074.

\end{thebibliography}

\end{document}